# A METHOD TO SOLVE THE DIOPHANTINE EQUATION $ax^2 - by^2 + c = 0$

Florentin Smarandache
University of New Mexico
200 College Road
Gallup, NM 87301, USA

**ABSTRACT**

We consider the equation
(1) $ax^2 - by^2 + c = 0$, with $a,b \in N^*$ and $c \in Z^*$.

It is a generalization of Pell's equation: $x^2 - Dy^2 = 1$. Here, we show that: if the equation has an integer solution and $a \cdot b$ is not a perfect square, then (1) has an infinitude of integer solutions; in this case we find a closed expression for $(x_n, y_n)$, the general positive integer solution, by an original method. More, we generalize it for any Diophantine equation of second degree and with two unknowns.

**INTRODUCTION**

If $ab = k^2$ is a perfect square ($k \in N$) the equation (1) has at most a finite number of integer solutions, because (1) become:
(2) $(ax - ky)(ax + ky) = -ac$.

If $(a,b)$ does not divide $c$, the Diophantine equation hasn't solutions.

**METHOD TO SOLVE**. Suppose (1) has many integer solutions.

Let $(x_0, y_0), (x_1, y_1)$ be the smallest positive integer solutions for (1), with $0 \leq x_0 < x_1$. We construct the recurrent sequences:

(3) $\begin{cases} x_{n+1} = \alpha x_n + \beta y_n \\ y_{n+1} = \gamma x_n + \delta y_n \end{cases}$

putting the condition (3) verify (1). It results:

$\begin{cases} a\alpha\beta = b\gamma\delta & (4) \\ a\alpha^2 - b\gamma^2 = a & (5) \\ a\beta^2 - b\delta^2 = -b & (6) \end{cases}$

having the unknowns $\alpha, \beta, \gamma, \delta$.

We pull out $a\alpha^2$ and $a\beta^2$ from (5), respectively (6), and replace them in (4) at the square; it obtains

$a\delta^2 - b\gamma^2 = a$. (7)

We subtract (7) from (5) and find

$\alpha = \pm\delta$. (8)

Replacing (8) in (4) it obtains

$\beta = \pm\dfrac{b}{a}\gamma$. (9)

Afterwards, replacing (8) in (5), and (9) in (6) it finds the same equation:

$a\alpha^2 - b\gamma^2 = a$. (10)

Because we work with positive solutions only, we take

$\begin{cases} x_{n+1} = \alpha_0 x_n + \dfrac{b}{a}\gamma_0 y_n \\ y_{n+1} = \gamma_0 x_n + \alpha_0 y_n \end{cases}$ ;

where $(\alpha_0, \gamma_0)$ is the smallest, positive integer solution of (10) such that $\alpha_0 \gamma_0 \neq 0$. Let $A = \begin{pmatrix} \alpha_0 & \dfrac{b}{a}\gamma_0 \\ \gamma_0 & \alpha_0 \end{pmatrix} \in M_2(Z)$.

Of course, if $(x', y')$ is an integer solution for (1), then $A\begin{pmatrix} x' \\ y' \end{pmatrix}$, $A^{-1}\begin{pmatrix} x' \\ y' \end{pmatrix}$ are another ones -- where $A^{-1}$ is the inverse matrix of $A$, i.e. $A^{-1} \cdot A = A \cdot A^{-1} = I$ (unit matrix). Hence, if (1) has an integer solution it has an infinite ones. (Clearly $A^{-1} \in M_2(Z)$).

The general positive integer solution of the equation (1) is
$(x'_n, y'_n) = (|x_n|, |y_n|)$.

$(GS_1)$ with $\begin{pmatrix} x_n \\ y_n \end{pmatrix} = A^n \cdot \begin{pmatrix} x_0 \\ y_0 \end{pmatrix}$, for all $n \in Z$,

where by conversion $A^0 = I$ and $A^{-k} = A^{-1} \cdots A^{-1}$ of $k$ times.

In problems it is better to write $(GS)$ as

$(GS_2)$ and
$\begin{pmatrix} x'_n \\ y'_n \end{pmatrix} = A^n \cdot \begin{pmatrix} x_0 \\ y_0 \end{pmatrix}$, $n \in N$

$\begin{pmatrix} x''_n \\ y''_n \end{pmatrix} = A^n \cdot \begin{pmatrix} x_1 \\ y_1 \end{pmatrix}$, $n \in N^*$

We proof, by *reduction ad absurdum*, $(GS_2)$ is a general positive integer solution for (1).

Let $(u, v)$ be a positive integer particular solution for (1). If

$$\exists k_0 \in N : (u,v) = A^{k_0} \begin{pmatrix} x_0 \\ y_0 \end{pmatrix}, \text{ or } \exists k_1 \in N^* : (u,v) = A^{k_1} \begin{pmatrix} x_1 \\ y_1 \end{pmatrix} \text{ then}$$

$(u,v) \in (GS_2)$. Contrary to this, we calculate $(u_{i+1}, v_{i+1}) = A^{-1} \begin{pmatrix} u_i \\ v_i \end{pmatrix}$

for $i = 0, 1, 2, \ldots$ where $u_0 = u, v_0 = v$. Clearly $u_{i+1} < u_i$ for all $i$. After a certain rank $x_0 < u_{i_0} < x_1$ it finds either $0 < u_{i_0} < x_0$ but that is absurd.

It is clear we can put

$$(GS_3) \begin{pmatrix} x_n \\ y_n \end{pmatrix} = A^n \cdot \begin{pmatrix} x_0 \\ \varepsilon y_0 \end{pmatrix}, \ n \in N, \text{ where } \varepsilon = \pm 1.$$

**We shall now transform the general solution $(GS_3)$ in a closed expression.**

Let $\lambda$ be a real number. $\text{Det}(A - \lambda \cdot I) = 0$ involves the solutions $\lambda_{1,2}$ and the proper vectors $V_{1,2}$ (i.e. $Av_i = \lambda_i v_i, \ i \in \{1,2\}$). Note $P = \begin{pmatrix} v_1 \\ v_2 \end{pmatrix}^t \in \mathcal{M}_2(\mathbf{R})$.

Then $P^{-1}AP = \begin{pmatrix} \lambda_1 & 0 \\ 0 & \lambda_2 \end{pmatrix}$, whence $A^n = P \begin{pmatrix} \lambda_1^n & 0 \\ 0 & \lambda_2^n \end{pmatrix} P^{-1}$, and replacing it in $(GS_3)$ and doing the calculus we find a closed expression for $(GS_3)$.

EXAMPLES

1. For the Diophantine equation $2x^2 - 3y^2 = 5$ at obtains

$$\binom{x_n}{y_n} = \begin{pmatrix} 5 & 6 \\ 4 & 5 \end{pmatrix}^n \cdot \binom{2}{3}, n \in N$$

and $\lambda_{1,2} = 5 \pm 2\sqrt{6},\ v_{1,2} = (\sqrt{6}, \pm 2)$, whence a closed expression for $x_n$ and $y_n$:

$$\begin{cases} x_n = \dfrac{4+\varepsilon\sqrt{6}}{4}(5+2\sqrt{6})^n + \dfrac{4-\varepsilon\sqrt{6}}{4}(5-2\sqrt{6})^n \\ y_n = \dfrac{3\varepsilon+2\sqrt{6}}{6}(5+2\sqrt{6})^n + \dfrac{3\varepsilon-2\sqrt{6}}{6}(5+-2\sqrt{6})^n \end{cases},$$

for all $n \in N$.

2. For equation $x^2 - 3y^2 - 4 = 0$ the general solution in positive integer is:

$$\begin{cases} x_n = (2+\sqrt{3})^n + (2-\sqrt{3})^n \\ y_n = \dfrac{1}{\sqrt{3}}\left[(2+\sqrt{3})^n + (2-\sqrt{3})^n\right] \end{cases},$$

for all $n \in N$, that is $(2,0), (4,2), (14,8), (52,30), \ldots$

EXERCICES FOR READER. Solve the Diophantine equations:

3. $x^2 - 12y^2 + 3 = 0$

$$\left[ \text{Remark: } \binom{x_n}{y_n} = \begin{pmatrix} 7 & 24 \\ 2 & 7 \end{pmatrix}^n \cdot \binom{3}{\varepsilon} = ?, n \in N \right]$$

4. $x^2 - 6y^2 - 10 = 0$.

$$\left[ \text{Remark:} \begin{pmatrix} x_n \\ y_n \end{pmatrix} = \begin{pmatrix} 5 & 12 \\ 12 & 5 \end{pmatrix}^n \cdot \begin{pmatrix} 4 \\ \varepsilon \end{pmatrix} = ?, n \in N \right]$$

5. $x^2 - 12y^2 - 9 = 0$

$$\left[ \text{Remark:} \begin{pmatrix} x_n \\ y_n \end{pmatrix} = \begin{pmatrix} 7 & 24 \\ 2 & 7 \end{pmatrix}^n \cdot \begin{pmatrix} 3 \\ 0 \end{pmatrix} = ?, n \in N \right]$$

6. $14x^2 - 3y^2 - 18 = 0$

## GENERALIZATIONS

If $f(x, y) = 0$ is a Diophantine equation of second degree and with two unknowns, by linear transformations it becomes

(12) $ax^2 + by^2 + c = 0$, with $a, b, c \in Z$.

If $ab \geq 0$ the equation has at most a finite number of integer solutions which can be found by attempts.

It is easier to present an example:

7. The Diophantine equation
(13) $9x^2 + 6xy - 13y^2 - 6x - 16y + 20 = 0$
can becomes
(14) $2u^2 - 7v^2 + 45 = 0$, where
(15) $u = 3x + y - 1$ and $v = 2y + 1$

We solve (14). Thus:

(16) $\begin{cases} u_{n+1} = 15u_n + 28v_n \\ v_{n+1} = 8u_n + 15v_n \end{cases}, n \in N$ with $(u_0, v_0) = (3, 3\varepsilon)$

**First solution:**

By induction we proof that: for all $n \in N$ we have $v_n$ is odd, and $u_n$ as well as $v_n$ are multiple of 3. Clearly $v_0 = 3\varepsilon, u_0$. For $n+1$ we have: $v_{n+1} = 8u_n + 15v_n = $ even + odd = odd, and of course $u_{n+1}, v_{n+1}$ are multiples of 3 because $u_n, v_n$ are multiple of 3, too.

Hence, there exist $x_n, y_n$ in positive integers for all $n \in N$:

(17) $\begin{cases} x_n = (2u_n - v_n + 3)/6 \\ y_n = (v_n - 1)/2 \end{cases}$

(from (15)). Now we find the $(GS_3)$ for (14) as closed expression, and by means of (17) it results the general integer solution of the equation (13).

**Second solution**

Another expression of the $(GS_3)$ for (13) we obtain if we transform (15) as: $u_n = 3x_n + y_n - 1$ and $v_n = 2y_n + 1$, for all $n \in N$. Whence, using (16) and doing the calculus, it finds

(18) $\begin{cases} x_{n+1} = 11x_n + \dfrac{52}{3}y_n + \dfrac{11}{3} \\ y_{n+1} = 12x_n + 19y_n + 3 \end{cases}$, $n \in N$,

with $(x_0, y_0) = (1,1)$ or $(2,-2)$ (two infinitude of integer solutions). Let

$A = \begin{pmatrix} 11 & 52/3 & 11/3 \\ 12 & 19 & 3 \\ 0 & 0 & 1 \end{pmatrix}$ Then $\begin{pmatrix} x_n \\ y_n \\ 1 \end{pmatrix} = A^n \begin{pmatrix} 1 \\ 1 \\ 1 \end{pmatrix}$

or

$$\begin{pmatrix} x_n \\ y_n \\ 1 \end{pmatrix} = A^n \begin{pmatrix} 2 \\ -2 \\ 1 \end{pmatrix}, \text{ always } n \in N. \qquad (19)$$

From (18) we have always $y_{n+1} \equiv y_n \equiv \cdots \equiv y_0 \equiv 1 \pmod{3}$, hence always $x_n \in Z$. Of course, (19) and (17) are equivalent as general integer solution for (13).

[The reader can calculate $A^n$ (by the same method liable to the start on this note) and find a closed expression for (19).]

**More generally:**

This method can be generalized for the Diophantine equations

(20) $\qquad \sum_{i=1}^{n} a_i X_i^2 = b, \text{ will all } a_i, b \text{ in } Z.$

If always $a_i a_j \geq 0, \ 1 \leq i \leq j < n$, the equation (20) has at most a finite number of integer solution.

Now, we suppose $\exists \ i_0, j_0 \in \{1, \ldots, n\}$ for which $a_{i_0} a_{j_0} < 0$ (the equation presents at least a variation of sign). Analogously, for $n \in N$. We define the recurrent sequences:

(21) $\qquad x_h^{(n+1)} = \sum_{i=1}^{n} a_{ih} x_i^{(n)}, \quad 1 \leq h \leq n$

considering $(x_1^0, \ldots, x_n^0)$ the smallest positive integer solution of (20). It replaces (21) in (20), it identifies the coefficients and it look for the $n^2$ unknowns $a_{ih}, \ 1 \leq i, h \leq n$. (This calculus is very intricate, but it can be done by means of a computer.) The method goes on similarly, but the calculus becomes more and more intricate - for example to calculate $A^n$. It must a computer, may be.

(The reader will be able to try his force for the Diophantine equation $ax^2 + by^2 - cz^2 + d = 0$, with $a, b, c \in N^*$ and $d \in Z$).

F. Smarandache, Sur la résolution d'équations du second degré a deux inconnues dans $Z$, in the book Généralizations et généralités, Ed. Nouvelle, Fès, Marocco; MR:85h:00003.

[Published in "Gaceta Matematica", 2a Serie, Volumen 1, Numero 2, 1988, pp.151-7; Madrid; translated in Spanish by Francisco Bellot Rasado: «Un metodo de resolucion de la ecuacion diofantica».]